\documentclass{article}

\usepackage{amssymb}
\usepackage{amsmath}
\usepackage{stmaryrd}
\usepackage{float}
\usepackage{graphicx}
\usepackage{tikz}
\usepackage{cite}

\usetikzlibrary{shapes.geometric, fit, positioning, arrows}
\usepackage{hyperref}
\usepackage{subcaption}

\newtheorem{assumption}{Assumption}

\usepackage[colors]{optsys}

\renewcommand{\va}[1]{{\boldsymbol{\uppercase{#1}}}}
\renewcommand{\opt}{{\sharp}}

\let\citep=\cite
\newcommand{\price}{p}

\newcommand{\alloc}{r}

\newcommand{\cF}{\mathcal{F}}

\newcommand{\st}{\text{s.t.}}

\newcommand{\dynamic}{{g}}
\newcommand{\x}{\boldsymbol{X}}
\renewcommand{\u}{\boldsymbol{U}}
\newcommand{\w}{\boldsymbol{W}}

\newcommand{\pvaluefunc}{\underline{V}}
\newcommand{\qvaluefunc}{\overline{V}}

\newcommand{\kk}{{(k)}}
\newcommand{\kp}{{(k+1)}}

\renewcommand{\horizon}{T}
\newcommand{\final}{T}
\renewcommand{\feedback}{\gamma}                              
\renewcommand{\policy}{\feedback}                             

\newcommand\produ{{\mathcal{N}}}
\newcommand\transport{{\mathcal{A}}}

\newcommand{\imag}{\mathrm{im}}

\newcommand{\pscal}[2]{#1 \cdot #2}
\renewcommand{\pscal}[2]{\big\langle#1\:,#2\big\rangle}     
\newcommand{\tstep}{\Delta T}

\newcommand{\card}[1]{\vert #1 \vert}
\newcommand{\NODES}{\mathcal{N}} 
\newcommand{\ARCS}{\mathcal{A}}

\newcommand{\VSDDP}{\underline{V}^{\mathrm{sddp}}}
\newcommand{\FORALLTIMES}[3]{\forall #1\in\ic{#2,#3}}

\newcommand{\post}{{t+1}}

\newcommand{\ARCFLOW}{Q}

\newcommand{\NODEFLOW}{F}
\newcommand{\IncidenceMatrix}{C}
\newcommand{\INCIDENCEMATRIX}{\mathcal{C}}

\def\graphicscale{0.48}
\graphicspath{{./}}
\usepackage{fullpage}
\title{Decentralized Multistage Optimization of Large-Scale Microgrids under Stochasticity}

\author{Fran\c{c}ois~Pacaud\thanks{MCS, Argonne National Laboratory, Lemont, USA}
  \and
  Michel~De~Lara\thanks{CERMICS, Ecole des Ponts, Marne-la-Vall\'ee, France}
  \and
  Jean-Philippe~Chancelier\footnotemark[2]
  \and
  Pierre~Carpentier\thanks{UMA, ENSTA Paris, IP Paris, France}
}

\date{\today}

\begin{document}
\renewcommand{\labelitemi}{\textbullet}

\maketitle
\begin{abstract}
  Microgrids are recognized as a relevant tool to absorb decentralized
  renewable energies in the energy mix. However, the sequential handling of
  multiple stochastic productions and demands, and of storage, make their
  management a delicate issue. We add another   layer of complexity by
  considering microgrids where different buildings stand   at the nodes of
  a network and are connected by the arcs; some buildings host local production
  and storage capabilities, and can exchange with others their energy surplus.
  We formulate the problem as a multistage stochastic optimization problem,
  corresponding to the minimization of the expected temporal sum of operational
  costs, while   satisfying the energy demand at each node, for all time.
  The resulting mathematical problem has a large-scale nature, exhibiting both
  spatial   and temporal couplings. However, the problem displays   a network
  structure that makes it   amenable to a mix of spatial
  decomposition-coordination with temporal decomposition methods.
  We conduct numerical simulations on microgrids of different sizes and
  topologies, with up to 48~nodes and 64~state variables.
  Decomposition methods are faster and provide more efficient policies
  than a state-of-the-art Stochastic Dual Dynamic Programming algorithm.
  Moreover, they scale almost linearly with the
  state dimension, making them a promising tool to address more complex
  microgrid optimal management problems.
\end{abstract}

\section{Introduction}

\subsection{Problem statement}

Power networks are organized more and more in a decentralized fashion,
with microgrids coordinating the production of local renewable energies
integrated with distributed storage. 
A broad overview of the emergence of consumer-centric
electricity markets is given in~\cite{sousa2019peer}, and challenges
associated with the integration of renewable energy can be found
in~\cite{morales2013integrating}.

As renewable energies, like sun and wind, are stochastic,
the Energy Management System (EMS) problem can naturally be formulated as
a multistage stochastic optimization problem~\cite{zheng2014stochastic}.
If we consider a microgrid consisting of ten buildings each equipped with a hot
water tank and a battery, controlled every quarter of an hour during one day,
such a problem is already large-scale (a hundred stages and a state with
dimension twenty) and a direct resolution is out of reach.
The large-scale nature of the EMS problem makes spatial decomposition methods
appealing~\cite{morstyn2018}.

Different distributed variants of the
Model Predictive Control algorithm have been proposed to control
microgrids~\cite{stadler2016distributed}.
Decomposition methods are also adapted to the resolution of
large-scale unit-commitment problems~\cite{bacaud2001bundle,papavasiliou2013multiarea}.
We refer the reader to~\cite{molzahn2017survey} for a recent survey of
distributed optimization methods for electric power systems, and to
\cite{kraning2014dynamic} for an example of a distributed optimization
algorithm applied to the control of a large power network.
When the system is dynamical, as is the case with storage, temporal
decomposition methods are also appealing.
Stochastic Dynamic Programming (SDP)~\cite{bellman57} is the reference algorithm but is
hampered by the well-known curse of dimensionality when the state dimension
exceeds five. Stochastic Dual Dynamic Programming (SDDP) takes the relay
under linear-convex assumptions with nice results,
for example in dam management~\citep{shapiro2012final}, but is efficient up to a certain
state dimension.

However, despite limitations inherent to SDP, one can go further
by mixing 
decomposition methods together.
For example, recent developments have mixed, in a stochastic setting, spatial decomposition methods with SDP
to effectively solve large-scale multistage stochastic optimization problems, by means
of the so-called Dual Approximate Dynamic Programming (DADP)
algorithm~\citep{barty2010decomposition}. 
In this paper, the EMS problems that we consider display a structure that
makes them amenable to a mix of spatial and temporal decomposition methods,
as developed previously in \cite{thesepacaud,Carpentier-Chancelier-DeLara-Pacaud:2020}
for general coupling constraints.
Indeed the (global) problem is naturally formulated
as a sum of local multistage stochastic optimization subproblems
coupled together via the global network constraints (flow
conservation on the graph).

\subsection{Contributions}

The contributions of this article are threefold.
i)
We implement DADP in an extended framework that supports generic
coupling constraints, specified on a directed graph
(whereas the previous implementations of DADP~\citep{carpentier2018stochastic}
considered only a unique central coupling constraints, or coupling
constraints formulated on a tree).
ii)
We implement a new algorithm called
\emph{Primal Approximate Dynamic Programming} (PADP), based on resource decomposition.
iii)
Thus equipped, 
on the one hand, we readily compute
an exact upper bound (PADP) and an exact lower bound (DADP)
of the global minimization problem and, on the other hand,
we yield two online control policies implementable by the EMS.
We provide numerical comparisons with a state-of-the-art SDDP algorithm,
and we show the effectiveness
of the two decomposition algorithms: for problems with more than 12 nodes,
the decomposition algorithms
converge faster than SDDP and yield control policies with lower costs.

\subsection{Structure of the paper}
The paper is organized as follows.
In Sect.~\ref{Management_of_large_scale_microgrids}, we outline the class of
optimal energy management problems that we address,
and in Sect.~\ref{Mathematical_formulation} we present
the associated mathematical formulation.
In Sect.~\ref{Resolution_by_decentralized_optimization},
we detail how to design algorithms by a mix of
spatial decomposition and of SDP.
In Sect.~\ref{chap:district:numerics}, we present the results of
numerical simulations for different microgrids of increasing size and complexity.

%
%

\section{Management of large-scale microgrids}
\label{Management_of_large_scale_microgrids}

Efficacity is the French urban Energy Transition Institute (ITE) devoted to
develop and implement innovative
solutions to build and manage energy-efficient cities.
It was created in 2014 with both the French government support and contributions
from companies, small and large. Efficacity has solicited us to
address the optimal energy management of urban electrical microgrids,
proposing several microgrid configurations,
both in terms of topology (structure of the districts) and in terms of equipment
(energy production and storage).

Buildings are heterogeneous: all are equipped with an electrical hot water tank,
but only some have solar panels and some others have batteries.
Indeed, as batteries and solar panels are expensive, they are shared out across
the network.
All units have the possibility to import and export energy
to and from the other buildings.
Moreover, we suppose that, if the local production is unable to fulfill the local demand,
even after exchanges between buildings, energy can be imported from an
external (regional or national) grid as a recourse.
Thus, each building is a decision unit able to locally consume,
produce, store, and also to exchange energy with other units and with the
external grid. The flows in the microgrid are impacted by uncertainties, both in demand
(e.g. electrical) and in production (e.g. solar panels).
We suppose that all actors are benevolent, allowing
a central planner (namely the EMS) to coordinate the local units.
The EMS aims at satisfying the balance between production and demand at each
node (building), but at least cost.

We manage the microgrids over one day, with decisions taken every 15 min.
The local solar energy productions match realistic data
corresponding to a summer day in Paris. The local demands
are generated using a stochastic simulator experimentally
validated~\cite{schutz2015comparison}.


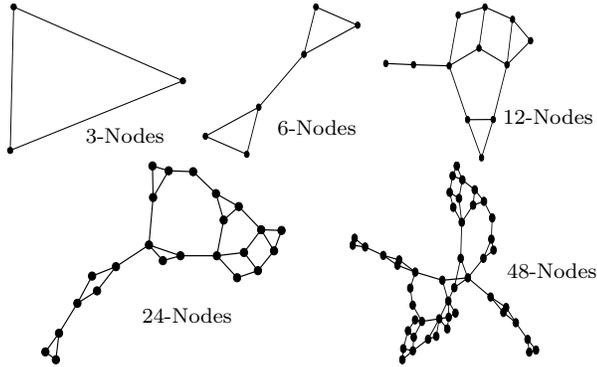
\begin{figure}[!ht]
  \begin{center}
\begin{picture}(0,0)%
\includegraphics{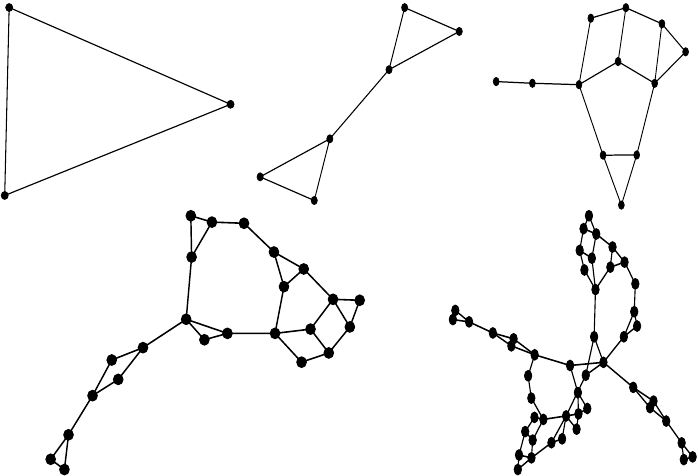}%
\end{picture}%
\setlength{\unitlength}{2072sp}%
\begingroup\makeatletter\ifx\SetFigFont\undefined%
\gdef\SetFigFont#1#2#3#4#5{%
  \reset@font\fontsize{#1}{#2pt}%
  \fontfamily{#3}\fontseries{#4}\fontshape{#5}%
  \selectfont}%
\fi\endgroup%
\begin{picture}(6373,4320)(901,-5011)
\put(1846,-2356){\makebox(0,0)[lb]{\smash{{\SetFigFont{8}{9.6}{\familydefault}{\mddefault}{\updefault}{\color[rgb]{0,0,0}$3$-Nodes}%
}}}}
\put(2521,-4516){\makebox(0,0)[lb]{\smash{{\SetFigFont{8}{9.6}{\familydefault}{\mddefault}{\updefault}{\color[rgb]{0,0,0}$24$-Nodes}%
}}}}
\put(6886,-3976){\makebox(0,0)[lb]{\smash{{\SetFigFont{8}{9.6}{\familydefault}{\mddefault}{\updefault}{\color[rgb]{0,0,0}$48$-Nodes}%
}}}}
\put(4141,-2266){\makebox(0,0)[lb]{\smash{{\SetFigFont{8}{9.6}{\familydefault}{\mddefault}{\updefault}{\color[rgb]{0,0,0}$6$-Nodes}%
}}}}
\put(6841,-2131){\makebox(0,0)[lb]{\smash{{\SetFigFont{8}{9.6}{\familydefault}{\mddefault}{\updefault}{\color[rgb]{0,0,0}$12$-Nodes}%
}}}}
\end{picture}%

    \caption{Examples of microgrid topologies\label{fig:nodal:graphtopology}}
  \end{center}
\end{figure}

We consider five different microgrids with growing sizes and different
topologies, that is, different nodes and connecting arcs.
The structure of the microgrids (see Figure~\ref{fig:nodal:graphtopology})
as well as the location of batteries and solar panels
come from case studies provided by Efficacity.
As an example (see Table~\ref{tab:numeric:pbsize}), the 12-Nodes problem consists of 12 buildings;
4 buildings are equipped with a 3~kWh battery, and 4 other
buildings are equipped with 16~$\text{m}^2$ of solar panels.
The devices are dispatched so that a building equipped with a solar
panel is connected to at least one building with a battery.

%

\section{Mathematical formulation}
\label{Mathematical_formulation}

We now address the mathematical formulation of a multistage stochastic
optimization problem that corresponds to
satisfying the supply-demand constraint at each node, at least expected cost.
After describing the basic mathematical objects
in~\S\ref{Description_of_the_problems},
we write mathematical equations related to arcs
in~\S\ref{Equations_related_to_arcs},
to nodes in~\S\ref{Equations_related_to_nodes},
and finally present the multistage stochastic optimization problem
formulation in~\S\ref{Formalization_as_a_multistage_stochastic_optimization_problem}.

\subsection{Network, stages and probability space}
\label{Description_of_the_problems}

We represent a district microgrid by a directed graph
$(\NODES, \ARCS)$, with $\NODES$ the set of nodes and $\ARCS$
the set of arcs. We denote by $\card{\NODES}$ the number of nodes,
and by $\card{\ARCS}$ the number of arcs.
We suppose that decisions are made at regular stages,
labeled by \( t \in \ic{0,\final} \),
where $\final \in \NN^\star$ is a finite horizon,
and where we use the notation \( \ic{r,s}=\na{r,r+1,\ldots,s-1,s} \) for two integers
$r \leq s$.
We have $\final = 96$ for a daily management
with decisions taken every 15 min.

We introduce a probability space $(\Omega, \cF, \PP)$,
denote the mathematical expectation by~$\EE$,
and write all random variables with uppercase bold letters.

\subsection{Equations related to arcs}
\label{Equations_related_to_arcs}


\subsubsection{Exchanging flows through arcs}
during the time interval~$[t,t+1)$,
each arc~$\arc \in \ARCS$ transports an energy flow~$\va\ARCFLOW^{\arc}_t\in\RR$, and each node
$\node \in \NODES$ imports or exports a flow $\va\NODEFLOW^{\node}_t\in\RR$.
The node flows~$\va\NODEFLOW^{\node}_t$ and the arc flows~$\va\ARCFLOW^{\arc}_t$
are related via a balance equation (Kirchhoff's current law)
written in matrix form as
\( \IncidenceMatrix \va\ARCFLOW_t + \va\NODEFLOW_t = 0 \),
where $\va\NODEFLOW_t = (\va\NODEFLOW^\node_t)_{\node \in \NODES} \in \RR^\NODES$
is the family of node flows at time~$t$,
$\va\ARCFLOW_t = (\va\ARCFLOW^\arc_t)_{\arc \in \ARCS} \in \RR^\ARCS$
is the family of arc flows at time~$t$ and where
$\IncidenceMatrix \in \na{-1,0,1}^{\NODES \times \ARCS}$
is the node-arc incidence matrix of the directed graph $(\NODES, \ARCS)$.
We identify $\va\NODEFLOW_t$ with a vector in $\RR^{\card{\NODES}}$
and $\va\ARCFLOW_t$ with a vector in $\RR^{\card{\ARCS}}$,
$\IncidenceMatrix$ being identified with a matrix with $\card{\NODES}$ rows
and $\card{\ARCS}$ columns.
We denote by
$\va\NODEFLOW = (\va\NODEFLOW_t)_{t \in \ic{0,\final-1}} \in \RR^{\final \cdot \card{\NODES}}$
the vector of node flows and by
$\va\ARCFLOW = (\va\ARCFLOW_t)_{t \in \ic{0,\final-1}} \in \RR^{\final \cdot \card{\ARCS}}$
the vector of arc flows.


\subsubsection{Transportation cost on arcs}
a quadratic cost~\( l_t^{\arc}(\va\ARCFLOW_t^{\arc}) =
c_2^{\arc} (\va\ARCFLOW_t^{\arc})^2 + c_1^{\arc} \va\ARCFLOW_t^{\arc} + c_0^{\arc}\)
(with given coefficients $c_2^\arc >0, c_1^\arc, c_0^\arc$,
so that each function~$l_t^\arc$ is strongly convex) is associated with
transporting the flow~$\va\ARCFLOW_t^{\arc}$ through arc~$\arc \in \ARCS$.
Such costs can arise from a difference in pricing, a fixed toll between
the different nodes, or by energy losses through the network.
We aggregate, in the \emph{global arc cost} 
\begin{equation}
  \label{eq:district:arccriterion}
  \Criterion_\transport(\va\ARCFLOW) =
  \EE \bgc{\sum_{\arc \in \ARCS} \sum_{t=0}^{\final-1}
  l_t^{\arc}(\va\ARCFLOW_t^{\arc})}
  \eqfinv
\end{equation}
all transport costs over all arcs in the graph,
over all times~$t\in\ic{0,\final-1}$
and over all random outcomes (hence the mathematical expectation term~$\EE$).

\subsection{Equations related to nodes}
\label{Equations_related_to_nodes}
We adopt a state space formalism to describe the physical equations related
to each node.

\subsubsection{State, control and uncertainty variables}
let $\na{\XX_t^{\node}}_{t\in \ic{0,\final}}$,
$\na{\UU_t^{\node}}_{t\in \ic{0,\final-1}}$ and
$\na{\WW_t^{\node}}_{t\in \ic{1,\final}}$
be sequences of Euclidean spaces of type~$\mathbb{R}^{p}$,
with appropriate dimensions~$p$ (possibly depending on time~$t$
and node~$\node \in \NODES$).

As all buildings hold a hot water tank, and some possibly also hold a battery,
the nodal state $\x_t^{\node}$ has dimension 1 or 2.
If the building at node~$\node$ hosts both a hot water tank
and a battery,
the state is $\x_t^{\node} = (\va b_t^{\node}, \va h_t^{\node})$ with values
in~$\XX_t^{\node}=\RR^2$,
where $\va b_t^{\node}$ (resp. $\va h_t^{\node}$) is the energy level inside the
battery (resp. hot water tank) at the beginning of the time
interval~$[t,t+1)$;
if the building at node~$\node$ only hosts a hot water tank,
the state is $\x_t^{\node} = \va h_t^{\node}$ with values in~$\XX_t^{\node}=\RR$.
The state at time~$0$ is supposed to be deterministic and known, equal to $x_0^{\node}$.

In the same way,
the nodal control $\u_t^{\node}$ has dimension 3 or 2.
If the building at node~$\node$ hosts a battery,
the control is $\u_t^{\node}=(\u_t^{b,\node},\u_t^{t,\node},\u_t^{ne,\node})$
with values in~$\UU_t^{\node}=\RR^3$,
where $\u_t^{b,\node}$ (resp. $\u_t^{t,\node}$) is
the amount of energy flowing into the
battery (resp. hot water tank), and $\u_t^{ne,\node}$ is the amount of electricity exchanged with
the external grid, during the time interval~$[t,t+1)$;
otherwise, the control is $\u_t^{\node}=(\u_t^{t,\node},\u_t^{ne,\node})$
with values in~$\UU_t^{\node}=\RR^2$.

Insofar, as the inhabitants of the different buildings have different
lifestyles, we suppose that each building has its own electrical and domestic
hot water demand profiles, and possibly its own solar panel production.
At node~$\node \in \NODES$, the uncertainty $\w_\post^{\node}=(\va{d}_\post^{hw,\node},\va{d}_\post^{el,\node})
\in \WW_\post^\node=\RR^2$ is made of the domestic hot water
demand~$\va{d}_\post^{hw,\node}$
and of~$\va{d}_\post^{el,\node}$, the local electricity demand minus the
production of the solar panel (if any),
both during the time interval~$[t,t+1)$.

\subsubsection{Dynamics inside each node}
we detail the dynamics in a building at node~$\node \in \NODES$.
We model the temporal evolution of a battery with the linear dynamics:
$\forall t \in \ic{0,\final-1}$,\footnote{%
    We have used the notation
$f^{+} = \max\na{0, f}$ and $f^{-} = \max\na{0, -f}$.}
\begin{subequations}
  \begin{equation}
    \label{eq:node:battery_dynamic}
    \va b_\post^{\node} = \alpha_b \va b_t^{\node} + \tstep
    \Bp{ \rho_c (\va u^{b,\node}_t)^+ - \dfrac{1}{\rho_d} (\va u^{b, \node}_t)^- }
    \eqfinv
  \end{equation}
  where $\alpha_b$ is the auto-discharge rate and $(\rho_d, \rho_c)$
  are given yields.
  We model the temporal evolution of an
  electrical hot water tank
  with the linear dynamics: $\forall t \in \ic{0,\final-1}$,
  \begin{equation}
    \label{eq:node:tank_dynamic}
    \va h_\post^{\node} = \alpha_h \va h_t^{\node} + \tstep
    \bp{ \beta_h \u_t^{t,\node} - \va d^{hw, \node}_\post }
    \eqfinv
  \end{equation}
  where 
  $\alpha_h$ is a discharge rate corresponding to losses by conduction
  and $\beta_h$ is a conversion coefficient.
\end{subequations}

We gather equations~\eqref{eq:node:battery_dynamic}-\eqref{eq:node:tank_dynamic}
in a \emph{nodal dynamics} function $\dynamic_t^\node: \XX_t^\node \times \UU_t^\node \times
\WW_\post^\node \rightarrow \XX_\post^\node$, stating that at each time $t$,
the next state $\x_\post^\node$ depends on the current state $\x_t^\node$, the
current decision $\u_t^\node$ and the uncertainty $\w_\post^\node$ occurring between
time $t$ and $t+1$.

\subsubsection{Load balance}
the load balance between production and demand
at node~$\node \in \NODES$ writes
\( \u_t^{ne,\node} - \va{d}_\post^{el,\node} - \u_t^{b,\node} - \u_t^{t,\node}
= \va\NODEFLOW_t^{\node} \),
where 
we recall that \( \va\NODEFLOW_t^{\node} \) is the energy exchanged with
the adjacent nodes.
Thus, for each time~$t  \in \ic{0,\final-1} $ and node $\node\in \NODES$, we introduce
the \emph{nodal load balance function}
\begin{equation}
  \label{eq:loadbalance}
  \Delta_t^{\node}(\x_t^{\node}, \u_t^{\node}, \w_\post^{\node}) =
  \u_t^{ne,\node} - \va{d}_\post^{el,\node} - \u_t^{b,\node} - \u_t^{t,\node}
  \eqfinp
\end{equation}
If $\Delta_t^\node(\cdot) <0$ (resp. $\Delta_t^\node(\cdot) >0$), the node
$\node$ imports (resp. exports) energy from (resp. to) adjacent nodes

\subsubsection{Cost function}
during the time interval~$[t,t+1)$ and at node $\node \in \NODES$,
the cost 
\begin{equation}
  L^{\node}_t(\va X_t^{\node},\va U_t^{\node},\va W^{\node}_{t+1})
  = p_t^{el} \u_t^{ne,\node}
\end{equation}
depends linearly on the price~$p_t^{el}$ to import electricity
from the external grid.
When we sum over time, we add a final penalization term~$K^{\node}(\va X_\final^{\node})$
to avoid an empty
electrical hot water tank at the end of the day.

\subsection{Multistage stochastic optimization problem formalization}
\label{Formalization_as_a_multistage_stochastic_optimization_problem}

We write
$\va\NODEFLOW^{\node} = (\va\NODEFLOW_0^{\node}, \cdots, \va\NODEFLOW^{\node}_{\final-1})\transp$
the node flow process arriving at each node $\node \in \NODES$ between times $0$ and $\final-1$.
We call optimal nodal cost the expression
\begin{subequations}
  \label{eq:district:localgenpb}
  \begin{align}
    \Criterion_\produ^{\node}
    & \np{\va\NODEFLOW^{\node},x_0^{\node}} = \nonumber \\
    & \min_{\x^{\node}, \u^{\node}} \;
      \EE \bgc{\sum_{t=0}^{\final-1}
      L^{\node}_t(\va X_t^{\node}, \va U_t^{\node}, \va W^{\node}_{t+1}) +
      K^{\node}(\x_\horizon^{\node})}
      \label{eq:nodalcostequation}
      \eqfinv
    \\
    & \quad \st\ \quad
      \FORALLTIMES{t}{0}{\final\!-\!1 }
      \nonumber
    \\
    & \phantom{\quad\st\ } \x_{t+1}^{\node} = \dynamic_t^{\node}(\x^{\node}_t, \u_t^{\node}, \w_{t+1}^{\node})
      \eqsepv \x_0^{\node} = x_0^{\node} \label{eq:dynamicsequation} \eqfinv
    \\
    & \phantom{\quad\st\ } \Delta_t^{\node}(\x_t^{\node}, \u_t^{\node}, \w_\post^{\node}) = \va\NODEFLOW_t^{\node}
    \label{eq:localbalanceequation} \eqfinv \\
    &\phantom{\quad\st\ } \sigma(\u_t^{\node}) \subset \sigma(\w_1, \cdots, \w_t, \w_\post)
    \label{eq:hazarddecisionframework} \eqfinv
  \end{align}
\end{subequations}
where $\x^{\node} = (\x^{\node}_0,\cdots,\x^{\node}_\final)$,
$\u^{\node}= (\u^{\node}_0,\cdots,\u^{\node}_{\final-1}) $
and $\w^{\node} = (\w^{\node}_0,\cdots,\w^{\node}_\final)$
are respectively the \emph{local} state (stocks), control (production)
and uncertainty (consumption) processes.
Similarly, we denote by~$\w = (\w^\node)_{\node\in\NODES}$
the \emph{global} uncertainty process.

To be able to almost surely satisfy the load balance
equations~\eqref{eq:localbalanceequation} at each node $\node \in \NODES$,
we assume that all decisions follow the \emph{hazard-decision}
information structure, that is, decision $\u_t^{\node}$ is taken
after the global uncertainty~$\w_{t+1} = (\w_\post^\node)_{\node\in\NODES}$
has been observed, hence the specific
form of the nonanticipativity constraint~\eqref{eq:hazarddecisionframework},
where $\sigma(\u_t^{\node})$ denotes the $\sigma$-algebra generated by the
random variable~$\u_t^{\node}$.The multistage nature of
Problem~\eqref{eq:district:localgenpb} stems from the nonanticipativity
constraint~\eqref{eq:hazarddecisionframework}.
Indeed, \eqref{eq:hazarddecisionframework} ensures that the decisions taken at time~$t$ depend
\emph{only} on the previous global uncertainties or, alternatively,
that the solution of~\eqref{eq:district:localgenpb}
is given as a sequence of \emph{policies} $\gamma_1^\node, \ldots, \gamma_{\final-1}^\node$,
such that, for all $t=1, \ldots, \final-1$, $\gamma_t^\node: \WW_1 \times \cdots \times\WW_{t+1} \to \UU_t$
and $\u_t^\node = \gamma_t^\node(\w_1, \ldots, \w_\post)$.


We have stated a global arc criterion in~\eqref{eq:district:arccriterion}
and local nodal criteria in~\eqref{eq:district:localgenpb},
both depending on node and arc flows coupled by Kirchhoff's current law
$\IncidenceMatrix \va\ARCFLOW_t + \va\NODEFLOW_t = 0$,
at each time~$t\in \ic{0,\final-1}$.
We rewrite these constraints globally as 
$\INCIDENCEMATRIX \va\ARCFLOW + \va\NODEFLOW = 0$
involving the global node flow and arc flow processes~$\va\ARCFLOW$ and~$\va\NODEFLOW$,
and where the matrix~$\INCIDENCEMATRIX \in
\RR^{\final \cdot \card{\NODES}} \times \RR^{\final \cdot \card{\ARCS}}$
is a block-diagonal matrix with matrix~$\IncidenceMatrix$ as diagonal element.

We set \( \XX_{0}=\prod_{\node \in \NODES} \XX_{0}^{\node} \) and,
for any $x_0=(x_0^\node)_{\node\in \NODES} \in \XX_{0}$,
we formulate the \emph{global} optimization problem of the central
manager (EMS) as
\begin{subequations}
  \label{transportproblem}
  \begin{align}
    V_0\opt(x_0) = \min_{ \va\NODEFLOW,\va\ARCFLOW} \;
    & \sum_{\node\in \NODES} \Criterion_\produ^{\node}\np{\va\NODEFLOW^{\node},x_0^{\node}} +
      \Criterion_\transport(\va\ARCFLOW) \\
    \st\
    & \INCIDENCEMATRIX \va\ARCFLOW + \va\NODEFLOW = 0
      \eqfinp
      \label{eq:couplingcons}
  \end{align}
\end{subequations}
The Problem~\eqref{transportproblem} couples $\card{\NODES} + 1$ independent
criteria through Constraint~\eqref{eq:couplingcons}. As the resulting
criterion is additive and Constraint~\eqref{eq:couplingcons}
is affine, Problem~\eqref{transportproblem} has a nice form
to use decomposition-coordination methods.



\section{Resolution by distributed optimization}
\label{Resolution_by_decentralized_optimization}

As just detailed in Sect.~\ref{Mathematical_formulation},
the global Problem~\eqref{transportproblem} encompasses a family
of local multistage stochastic optimization subproblems,
coupled together via a transportation problem corresponding
to the flows exchanged through the graph. We now detail how
to solve~\eqref{transportproblem} in a distributed fashion.

In \S\ref{sec:district:nodal}, we decouple~\eqref{transportproblem}
node by node using either price or resource decomposition schemes.
In \S\ref{sec:district:algorithm}, we show how to find the
most appropriate deterministic price and resource processes.
Thus, we obtain two algorithms, each of them yielding nodal value functions
and, from these latter,
upper and lower bounds for the optimal cost and online control policies.

\subsection{Mixing nodal and time decomposition} 
\label{sec:district:nodal}

In~\cite{Carpentier-Chancelier-DeLara-Pacaud:2020}, we introduced
a generic framework to bound a global problem by decomposing
it into smaller local subproblems, easier to solve.
Problem~\eqref{transportproblem} lies in the generic framework
introduced in~\cite{Carpentier-Chancelier-DeLara-Pacaud:2020},
and the coupling equation
$\INCIDENCEMATRIX \va\ARCFLOW + \va\NODEFLOW = 0$ is a special
case of the generic coupling constraint of this framework.
Thus, to solve Problem~\eqref{transportproblem}, 
we first apply \emph{spatial decoupling}
into nodal and arc subproblems, and then apply the
\emph{temporal decomposition} induced by Dynamic Programming.

\subsubsection{Price decomposition of the global problem}
we follow the procedure introduced in
\cite[\S2.2]{Carpentier-Chancelier-DeLara-Pacaud:2020}
to solve Problem~\eqref{transportproblem} by price decomposition
and to provide a lower bound of its optimal value $V_0\opt(x_0)$.
In the case under study,
price decomposition follows from the dualization of 
Constraint~\eqref{eq:couplingcons} using a deterministic price coordination
process~$\price=\np{\price^\node}_{\node \in \NODES} \in \RR^{\final\cdot\card{\NODES}}$
as multiplier.
We define the \emph{global price value function}\footnote{%
In the expression \( \pvaluefunc\nc{\price}\np{x_0} \),
we use brackets~$\nc{\price}$ to indicate a parametric dependence,
whereas we use parenthesis~$(x_0)$ to indicate the argument of the
function~\( \pvaluefunc\nc{\price} : \XX_{0} \to \RR \).
We also use the notation \( \pvaluefunc\nc{\price}\np{\cdot} \)
to designate this function.}
\( \pvaluefunc\nc{\price} : \XX_{0} \to \RR \)
associated with
Problem~\eqref{transportproblem} by the following expression,
for all $x_0=(x_0^\node)_{\node\in \NODES} \in \XX_{0}$,
\begin{eqnarray}{rl}
  \label{eq:district:dualvf}
  \pvaluefunc\nc{\price}\np{x_0} =
  & \min_{\va\NODEFLOW, \va\ARCFLOW} \; \sum_{\node \in \NODES}
  \Criterion_\produ^{\node}\np{\va\NODEFLOW^{\node},x_0^{\node}} +
  \Criterion_\transport(\va\ARCFLOW) \nonumber \\
  & \hfill +
  \EE \bc{\pscal{\price}{\INCIDENCEMATRIX \va\ARCFLOW + \va\NODEFLOW}} \eqfinp
\end{eqnarray}
We observe in a straightforward manner that
  \begin{equation}
    \label{eq:district:localpricepb}
    \pvaluefunc\nc{\price}\np{x_0} =
    \sum_{\node \in \NODES}\pvaluefunc_\produ^{\node}\nc{\price^{\node}}\np{x_0^{\node}}
    + \pvaluefunc_\transport\nc{\price}
    \eqfinv
  \end{equation}
  that is, the global price value function~$\pvaluefunc\nc{\price}\np{\cdot}$ naturally
  decomposes into a family of \emph{nodal price value functions}
  \( \pvaluefunc_\produ^{\node}\nc{\price^{\node}}: \XX_0^{\node} \to\RR \),
  $\forall n \in \NODES$, given by
  \begin{subequations}
  \begin{equation}
    \label{eq:district:localnodalpricepb}
    \pvaluefunc_\produ^{\node}\nc{\price^{\node}}\np{x_0^{\node}} = \min_{\va\NODEFLOW^{\node}}
    \Criterion_\produ^{\node}\np{\va\NODEFLOW^{\node},x_0^{\node}} + \EE \bc{ \pscal{\price^{\node}}{\va\NODEFLOW^{\node}}}
    \eqfinv
  \end{equation}
  and an \emph{arc price value function}\footnote{%
    Which, to the difference
    of the nodal price value function~\eqref{eq:district:localnodalpricepb},
    does not depend on the initial state~$x_0$.}
  given by 
  \begin{equation}
    \label{eq:district:localarcpricepb}
    \pvaluefunc_\transport\nc{\price} = \min_{\va\ARCFLOW}\;
    \Criterion_\transport(\va\ARCFLOW) + \EE \bc{ \pscal{\INCIDENCEMATRIX^{\top}\price}{\va\ARCFLOW}} \eqfinp
  \end{equation}
\end{subequations}
What is more,  for all $\node \in \NODES$, the optimal value~$\pvaluefunc_\produ^{\node}\nc{\price^{\node}}(x_0^{\node})$
can be computed by Dynamic Programming under the so-called
\emph{white noise assumption}.
\begin{assumption}
  \label{hyp:independent}
  The global uncertainty process $\np{\w_1,\cdots,\w_{\final}}$
  consists of stagewise independent random variables.
\end{assumption}

For all node~$\node \in \NODES$ and for any price
$\price^{\node} \in\RR^{\final}$,
taking into account the expression~\eqref{eq:district:localgenpb}
of the nodal cost~$\Criterion_\produ^{\node}$,
we introduce the sequence
$\na{\pvaluefunc_{\produ,t}^{\node}\nc{\price^{\node}}\np{\cdot}}_{t\in\ic{0,\final}}$
of local price value functions defined by
$\pvaluefunc_{\produ, \final}^{\node}\nc{\price^{\node}}\np{\cdot} =
K^{\node}\np{\cdot}$ (final cost) and then, inductively,
for all $t \in \ic{0,\final}$ and $x_t^{\node} \in \XX_t^{\node}$, by

\begin{align*}
  \pvaluefunc_{\produ, t}^{\node}\nc{\price^{\node}}(x_t^{\node}) =
  & \min_{\x^{\node}, \u^{\node}, \va\NODEFLOW^{\node}}
    \EE \bigg[ \sum_{s=t}^{\final-1}
    \Big(L^{\node}_s(\va X_s^{\node}, \va U_s^{\node}, \va W^{\node}_{s+1})
  \\
  & \phantom{\min_{\x^{\node}, \u^{\node}, \va\NODEFLOW^{\node}}}
    + \pscal{\price_s^{\node}}{\va\NODEFLOW_s^{\node}} \Big)
  + K^{\node}(\x_\horizon^{\node})\bigg] \eqfinv
  \\
  & \st \;\;
    \va \x_t^{\node} = x_t^{\node} \eqsepv
  \eqref{eq:dynamicsequation}-\eqref{eq:localbalanceequation}-\eqref{eq:hazarddecisionframework}
  \eqfinp
\end{align*}
Under Assumption~\ref{hyp:independent}, these local price
value functions satisfy the following Dynamic Programming equations:
for all $\node \in \NODES$,
\(
    \pvaluefunc_{\produ, \final}^{\node}\nc{\price^{\node}}(x_\final^{\node})
    = K^{\node}(x_\final^{\node}) \eqfinv
    \)
  and, for \( t=\final\!-\!1, \ldots, 0 \),
  \begin{align}
    \label{eq:bellmanequationpricefunctions}
    \pvaluefunc_{\produ, t}^{\node}\nc{\price^{\node}}(x_t^{\node}) =
    & \EE
      \Big[\min_{u_t^{\node}}
      L^{\node}_t(x_t^{\node},u_t^{\node},\w^{\node}_{t+1}) \nonumber \\
    & + \pscal{\price_t^{\node}}
       {\Delta_t^{\node}(x_t^{\node},u_t^{\node},\w_{t+1}^{\node})} \\
    & + \pvaluefunc_{\produ, \post}^{\node}\nc{\price^{\node}}
       \bp{\dynamic_t^{\node}(x_t^{\node},u_t^{\node}, \w_{t+1}^{\node})}
       \Big]
       \eqfinp \nonumber
  \end{align}
The nodal price value function $\pvaluefunc_\produ^{\node}\nc{\price^{\node}}\np{\cdot}$
in~\eqref{eq:district:localnodalpricepb}
is equal to the local price value function at time~$t=0$:
\(
\pvaluefunc_\produ^{\node}\nc{\price^{\node}}\np{x_0^{\node}} =
\pvaluefunc_{\produ, 0}^{\node}\nc{\price^{\node}}(x_0^{\node}) 
\), for all for $x_0^{\node} \in \XX_0^{\node}$.


Considering the expression~\eqref{eq:district:arccriterion}
of the arc cost~$\Criterion_\transport(\va\ARCFLOW)$,
the arc price value function~$\pvaluefunc_\transport\nc{\price}$
is additive \wrt\ (with respect to) time and space, and thus can be decomposed
at each time~$t$ and each arc~$\arc$. The resulting arc
subproblems do not involve any time coupling and can be computed
by standard mathematical programming tools or even analytically.

\subsubsection{Resource decomposition of the global problem}

we now solve Problem~\eqref{transportproblem} by resource decomposition
(see \cite[\S2.2]{Carpentier-Chancelier-DeLara-Pacaud:2020}) using
a deterministic resource process
$\alloc=\np{\alloc^\node}_{\node \in \NODES} \in \RR^{\final\cdot\card{\NODES}}$,
such that $\alloc \in \imag(\INCIDENCEMATRIX)$.\footnote{If
  $\alloc\notin\imag(\INCIDENCEMATRIX)$, we have~$\qvaluefunc\nc{\alloc}=+\infty$ in~\eqref{eq:district:quantdec} as the constraint
  $\INCIDENCEMATRIX \va\ARCFLOW + \alloc = 0$ cannot be satisfied.}
We decompose the global constraint~\eqref{eq:couplingcons}
\wrt\ nodes and arcs as
$\va\NODEFLOW = \alloc , \INCIDENCEMATRIX \va\ARCFLOW = -\alloc$.
We define the \emph{global resource value function}~$\qvaluefunc\nc{\alloc}\np{\cdot}$ associated with
Problem~\eqref{transportproblem} by the following expression, for all
$x_0=(x_0^\node)_{\node\in \NODES} \in \XX_{0}$:
\begin{subequations}\label{eq:district:quantdec}
  \begin{align}
    \qvaluefunc\nc{\alloc}\np{x_0} =
    & \min_{\va\NODEFLOW, \va\ARCFLOW} \; \sum_{\node \in \NODES}
    \Criterion_\produ^{\node}\np{\va\NODEFLOW^{\node},x_0^{\node}} +
    \Criterion_\transport(\va\ARCFLOW) \eqsepv \\
    & \st\ \;\;
      \va\NODEFLOW - \alloc = 0 \eqsepv
    \INCIDENCEMATRIX \va\ARCFLOW + \alloc = 0 \eqfinp
  \end{align}
\end{subequations}
We observe in a straightforward manner that
  \begin{equation}
    \qvaluefunc\nc{\alloc}\np{x_0} =
    \sum_{\node \in \NODES}\qvaluefunc_\produ^{\node}\nc{\alloc^{\node}}\np{x_0^{\node}}
    +   \qvaluefunc_\transport\nc{\alloc} \eqfinv
  \end{equation}
  that is, the global resource value function~$\qvaluefunc\nc{\alloc}\np{\cdot}$ naturally
  decomposes into a family of \emph{nodal resource value functions}~\(
\qvaluefunc_\produ^{\node}\nc{\alloc^{\node}}\np{\cdot} \) defined,
for all $\node \in \NODES$ and $x_0^{\node} \in \XX_0^{\node}$, by
\begin{subequations}
  \begin{equation}
    \label{eq:district:localnodalallocpb}
    \qvaluefunc_\produ^{\node}\nc{\alloc^{\node}}\np{x_0^{\node}} = \min_{\va\NODEFLOW^{\node}} \; \Criterion_\produ^{\node}\np{\va\NODEFLOW^{\node},x_0^{\node}}
    \quad \st \quad \va\NODEFLOW^{\node} - \alloc^{\node} = 0
  \end{equation}
  and an \emph{arc resource value function}
  (not depending on~$x_0$)
  \begin{equation}
    \label{eq:district:localarcallocpb}
    \qvaluefunc_\transport\nc{\alloc}
    = \min_{\va\ARCFLOW} \; \Criterion_\transport(\va\ARCFLOW)
    \quad \st \quad \INCIDENCEMATRIX \va\ARCFLOW + \alloc = 0 \eqfinp
  \end{equation}
\end{subequations}
For all node~$\node \in \NODES$,
taking into account the expression~\eqref{eq:district:localgenpb}
of the nodal cost~$\Criterion_\produ^{\node}$,
we introduce the sequence
$\na{\qvaluefunc_{\produ,t}^{\node}\nc{\alloc^{\node}}\np{\cdot}}_{t\in\ic{0,\final}}$
of local resource value functions defined by
$\qvaluefunc_{\produ, \final}^{\node}\nc{\alloc^{\node}}\np{\cdot} =
K^{\node}\np{\cdot}$ (final cost) and then, inductively,
for all $t \in \ic{0,\final}$ and $x_t^{\node} \in \XX_t^{\node}$, by
\begin{align*} 
  \qvaluefunc_{\produ,t}^{\node}\nc{\alloc^{\node}}(x_t^{\node}) =
  & \min_{\substack{\x^{\node}\\ \u^{\node}\\ \va\NODEFLOW^{\node}}}
    \EE
    \bigg[ \sum_{s=t}^{\final-1}
    L^{\node}_s(\va X_s^{\node}, \va U_s^{\node}, \va W^{\node}_{s+1})
      {+}K^{\node}(\x_\horizon^{\node})\bigg]
      \nonumber 
  \\
  & \st\ \; \va x_t^{\node} = x_t^{\node} \eqsepv
  \eqref{eq:dynamicsequation}-\eqref{eq:localbalanceequation}-\eqref{eq:hazarddecisionframework}
  \eqfinv \\
  & \phantom{\st\ \;} \va\NODEFLOW^{\node}_s - \alloc_s^{\node} = 0 \eqsepv
  \FORALLTIMES{s}{t}{\final\!-\!1} \eqfinp
\end{align*}
If Assumption~\ref{hyp:independent} holds true,
$\qvaluefunc_\produ^{\node}\nc{\alloc^{\node}}(x_0^{\node})$
can be computed by Dynamic Programming and the nodal resource
value function $\qvaluefunc_\produ^{\node}\nc{\alloc^{\node}}\np{\cdot}$ in~\eqref{eq:district:localnodalallocpb}
is equal to the local resource value function at time~$t=0$:
\(
\qvaluefunc_\produ^{\node}\nc{\alloc^{\node}}\np{x_0^{\node}} =
\qvaluefunc_{\produ, 0}^{\node}\nc{\alloc^{\node}}(x_0^{\node})
\)   for all $x_0^{\node} \in \XX_0^{\node}$.
%
In the case of resource decomposition, arcs are coupled
through the constraint $\INCIDENCEMATRIX \va\ARCFLOW + \alloc = 0$, so that
the arc resource value function $\qvaluefunc_\transport\nc{\alloc}$
in \eqref{eq:district:localarcallocpb} is not additive
in space, but remain additive \wrt\ time. As in price
decomposition, it can be computed by standard mathematical
programming tools or even analytically.

\subsubsection{Upper and lower bounds of the global problem}
applying \cite[Proposition~2]{Carpentier-Chancelier-DeLara-Pacaud:2020}
to the global price value function~\eqref{eq:district:dualvf}
and resource value functions~\eqref{eq:district:quantdec},
we are able to bound up and down the optimal value $V_0\opt(x_0)$
of Problem~\eqref{transportproblem} 
as follows:
\begin{equation}
  \label{eq:district:upperlowerboundglobal}
  \pvaluefunc\nc{\price}\np{x_0}
  \leq V_0\opt(x_0) \leq
  \qvaluefunc\nc{\alloc}\np{x_0}
  \eqsepv  \forall x_0\in \prod_{\node \in \NODES} \XX_{0}^{\node}
  \eqfinp
\end{equation}
These inequalities hold true for any price process
$\price \in \RR^{\horizon \cdot \card{\NODES}}$
and for any resource process $\alloc \in \imag(\INCIDENCEMATRIX)
\subset \RR^{\horizon \cdot \card{\NODES}}$.

\subsection{Algorithmic implementation}
\label{sec:district:algorithm}

In \S\ref{sec:district:nodal}, we have decomposed
Problem~\eqref{transportproblem} spatially and temporally:
the global problem is now split into (small) subproblems using
price and resource decompositions, and each subproblem is
solved by Dynamic Programming. These decompositions yield
bounds for the value of the global problem. To obtain
tighter bounds for the optimal value in~\eqref{eq:district:upperlowerboundglobal},
we follow the approach presented
in~\cite[\S3.2]{Carpentier-Chancelier-DeLara-Pacaud:2020},
that is, we maximize (resp. minimize) the left-hand side
(resp. the right-hand side)
in \eqref{eq:district:upperlowerboundglobal} \wrt\
the price vector~$\price\in\RR^{\horizon \cdot \card{\NODES}}$
(resp. the resource vector~$\alloc\in\RR^{\horizon \cdot \card{\NODES}}$),  using a gradient-like algorithm.

\subsubsection{Lower bound improvement}
\label{subsec:district:algorithm:price}
we detail how to improve the lower bound given
by the price value function
in~\eqref{eq:district:upperlowerboundglobal}.
We fix
$x_0=(x_0^\node)_{\node\in \NODES} \in \XX_{0}$,
and we aim at solving
\(
    \max_{\price\in \RR^{\final\cdot\card{\NODES}}} \;
    \pvaluefunc\nc{\price}\np{x_0}
\),
that is, written equivalently (see \eqref{eq:district:dualvf})
\begin{equation}
  \label{eq:district:relaxedmpts}
  \max_{\price\in \RR^{\final\cdot\card{\NODES}}} \;
  \min_{\va\NODEFLOW, \va\ARCFLOW} \; \sum_{\node\in\NODES}
    \Criterion_\produ^{\node}\np{\va\NODEFLOW^{\node},x_0^{\node}} +
  \Criterion_\transport(\va\ARCFLOW)
    \hfill
  + \pscal{\price}{\EE\bc{\INCIDENCEMATRIX \va\ARCFLOW + \va\NODEFLOW}}
  \eqfinp
\end{equation}
We solve the maximization Problem~\eqref{eq:district:relaxedmpts}
\wrt\ $\price$ using a gradient ascent method (Uzawa algorithm).
At iteration~$k$, we suppose given a deterministic price process
$\price^\kk$ and a gradient step $\rho^\kk$.
The algorithm proceeds as follows,
$ \forall \node \in \NODES$,
\begin{subequations}
  \label{eq:district:subgradientpricedecom}
  \begin{align}
    \label{eq:flowupdate}
    {\va\NODEFLOW^{\node}}^\kp
    & \in \argmin_{\va\NODEFLOW^{\node}}
    \Criterion_\produ^{\node}\np{\va\NODEFLOW^{\node},x_0^{\node}} + \EE \bc{\pscal{{\price^{\node}}^\kk}{\va\NODEFLOW^{\node}}} \eqfinv \\
    \label{eq:arcupdate}
    \va\ARCFLOW^\kp
    & \in \argmin_{\va\ARCFLOW}
    \Criterion_\transport(\va\ARCFLOW) + \EE \bc{ \pscal{\INCIDENCEMATRIX^{\top}\price^\kk}{\va\ARCFLOW}}
    \eqfinv
    \\
    \price^\kp
    & = \price^\kk + \rho^\kk \; \EE\bc{\INCIDENCEMATRIX \va\ARCFLOW^\kp + \va\NODEFLOW^\kp} \eqfinp
    \label{eq:priceupdate}
  \end{align}
\end{subequations}
At each iteration $k$, updating $\price^\kk$ requires
the computation of the gradient
of~$\nabla_\price \pvaluefunc\nc{\price^\kk}\np{x_0}$,
that is, the expected value
$\EE\bc{\INCIDENCEMATRIX \va\ARCFLOW^\kp + \va\NODEFLOW^\kp}$,
usually estimated by a Monte-Carlo method.
The price update formula~\eqref{eq:priceupdate} --- corresponding
to the standard gradient algorithm for the maximization
\wrt\ $\price$ in Problem~\eqref{eq:district:relaxedmpts} ---
can be replaced by more sophisticated algorithms (e.g. quasi-Newton).

\subsubsection{Upper bound improvement}
\label{subsec:district:algorithm:alloc}
we now focus on the improvement of the upper bound given
by the global resource value function
in~\eqref{eq:district:upperlowerboundglobal}.
We fix
$x_0=(x_0^\node)_{\node\in \NODES} \in \prod_{\node \in \NODES} \XX_{0}^{\node}$,
and we aim at solving
\(
    \min_{\alloc \in \imag(\INCIDENCEMATRIX)}  \;
    \qvaluefunc\nc{\alloc}\np{x_0}
  \),
  whose detailed expression is
  (see \eqref{eq:district:quantdec}):
  \begin{align}
    \label{eq:district:quantdecompositionrelaxed}
    \min_{\alloc \in \imag(\INCIDENCEMATRIX)}
    \bigg(
    & \sum_{\node \in \NODES}
      \Bp{\min_{\va\NODEFLOW^{\node}} \; \Criterion_\produ^{\node}\np{\va\NODEFLOW^{\node},x_0^{\node}}
    \quad \st \quad \va\NODEFLOW^{\node} - \alloc^{\node} \!=\! 0}
    \nonumber \\
    &+ \Bp{\min_{\va\ARCFLOW} \; \Criterion_\transport(\va\ARCFLOW)
    \quad \st \quad \INCIDENCEMATRIX \va\ARCFLOW + \alloc \!=\! 0}\bigg) \eqfinp
  \end{align}
We now sketch how we solve the minimization problem~\eqref{eq:district:quantdecompositionrelaxed}
using a gradient-like method.
As we recognize in~\eqref{eq:district:quantdecompositionrelaxed}
the minimization problems~\eqref{eq:district:localnodalallocpb}
and~\eqref{eq:district:localarcallocpb},
we obtain that the gradients \wrt\ $\alloc$
in~\eqref{eq:district:quantdecompositionrelaxed}
can be expressed as mathematical expectations
$\mu^\node= \EE[\va{M}^\node]$
and $\xi = \EE[\va{\Xi}]$, where we have denoted
by~${\va{M}^\node}$ the optimal multiplier associated
with the constraint $\va\NODEFLOW^\node - \alloc^{\node} = 0$
in~\eqref{eq:district:localnodalallocpb}, and by~$\va{\Xi}$ the
optimal multiplier associated with the constraint $\INCIDENCEMATRIX \va\ARCFLOW  + \alloc = 0$
in~\eqref{eq:district:localarcallocpb}.

The algorithm proceeds as follows.
At each iteration~$k$, the algorithm updates the
resource $\alloc^\kk$ and the gradient step~$\rho^\kk$.
We solve the optimization problems~\eqref{eq:district:localnodalallocpb}
and~\eqref{eq:district:localarcallocpb} with the resource set to $\alloc^\kk$ and obtain
the optimal solutions together with the associated multipliers $\mu^{\node,\kk}$ and $\xi^{\kk}$
as described in the previous paragraph. Then, the resource $\alloc^\kp$ is updated by
\begin{align}
  \alloc^\kp
  & = \projop{\imag(\INCIDENCEMATRIX)} \Bp{\alloc^\kk - \rho^\kk \;
    \bp{\mu^\kp + \xi^\kp}} \eqfinv
  \label{eq:district:subgradientquantdecom}
\end{align}
where $\projop{\imag(\INCIDENCEMATRIX)}$ is the orthogonal projection
onto the subspace $\imag(\INCIDENCEMATRIX)$
and $\mu^\kp = \{\mu^{n,\kp}\}_{\node \in \NODES}$. The multipliers $\mu^n$ and $\xi$
are approximated using a Monte Carlo method.
Again, the above projected gradient algorithm, 
used to update the resource, can be replaced by any gradient-based constrained optimization algorithm.

\subsubsection{Upper and lower nodal value functions}

The two algorithms \eqref{eq:district:subgradientpricedecom} and
\eqref{eq:district:subgradientquantdecom}
converge respectively to a price process $\price^{\infty}$
and to a resource process $\alloc^{\infty} \in \imag{(\INCIDENCEMATRIX)}$,
parameterized by a \emph{fixed} initial state. However, by
applying the two inequalities~\eqref{eq:district:upperlowerboundglobal}
with both \emph{admissible} processes~$\price^\infty$ and $\alloc^\infty$,
we are able to bound the optimal value function~$V_0\opt(\cdot)$ globally:
$\pvaluefunc\nc{\price^\infty}(x_0) \leq V_0\opt(x_0) \leq \qvaluefunc\nc{\alloc^\infty}(x_0)$
for all $x_0 \in \XX_{0}$.


\section{Numerical simulation results}
\label{chap:district:numerics}
\begin{table}[!ht]
  \centering
  \resizebox{.38\textwidth}{!}{
    \begin{tabular}{|c|ccccc|}
      \hline
      Problem           & $\card{\NODES}$  & $\card{\ARCS}$  & $dim(\XX_t)$ & $dim(\WW_t)$ & $supp(\w_t)$ \\
      \hline
      \hline
      \textrm{3-Nodes}  & 3               & 3               & 4            & 6            & $10^3$    \\
      \textrm{6-Nodes}  & 6               & 7               & 8            & 12           & $10^6$    \\
      \textrm{12-Nodes} & 12              & 16              & 16           & 24           & $10^{12}$ \\
      \textrm{24-Nodes} & 24              & 33              & 32           & 48           & $10^{24}$ \\
      \textrm{48-Nodes} & 48              & 69              & 64           & 96           & $10^{48}$ \\
      \hline
    \end{tabular}
  }
  \caption{Microgrid management problems (cardinals and dimensions)\label{tab:numeric:pbsize}}
\end{table}
We will term \emph{Dual Approximate Dynamic Programming} (DADP)
the price decomposition algorithm described
in \S\ref{subsec:district:algorithm:price} and
\emph{Primal Approximate Dynamic Programming} (PADP) the resource
decomposition algorithm described in \S\ref{subsec:district:algorithm:alloc}.
We compare DADP and PADP with the well-known Stochastic Dual Dynamic
Programming (SDDP) algorithm (see~\cite{girardeau2014convergence}
and references inside) on the results that they yield on
five microgrid optimal management problems with growing sizes:
  Table~\ref{tab:numeric:pbsize} displays the features (different sizes
and dimensions) of the cases we consider for numerical simulations.
For this purpose, we detail the offline computation of value functions
in~\S\ref{ssec:nodalalgorithms}, the online control policies
in~\S\ref{ssec:nodalpolicies}, and we finally summarize the numerical simulation
results and compare the three algorithms in~\S\ref{Numerical_results}.

For the uncertainties, we generate scenarios for the demands at
each node using the generator presented in~\cite{baetens2016modelling}, from
which we add the production of the solar panel for the corresponding buildings.
Then, we model the process~$\np{\va\w_t^\node}_{t, \node}$ as nonstationary random variables that are
stagewise and node by node
probabilistically independent
with a finite probability distribution on the set~$\WW_t^\node=\RR^2$.
For each case, we consider a single initial
state~$x_0$, corresponding to standard configurations of the storages (like
minimal energy in a battery).

\subsection{Computing offline value functions}
\label{ssec:nodalalgorithms}

\subsubsection{With the SDDP algorithm}
\label{SDDP_on_the_global_problem}
to compute value functions, the SDDP algorithm is
not implementable in a straightforward manner. Indeed,
even if the cardinality of the support of each local random variable~$\w_t^{\node}$ remains low,
the cardinal of the support of the global uncertainty~$\w_t$ becomes huge as
the number~$\card{\NODES}$ of nodes grows (see Table~\ref{tab:numeric:pbsize}),
so that the exact computation of expectations, as required at each time
step during the backward pass of the SDDP algorithm (see~\cite{shapiro10}),
becomes untractable.
To overcome this issue, we resample the probability distribution
of~$\np{\w_t^\node}_{\node \in \NODES}$ for each time~$t$,
to deal with an uncertainty support of reasonable size, using the
$k$-means clustering method in~\cite{rujeerapaiboon2018scenario}.
As the local problems are convex \wrt\ the uncertainties, by
Jensen inequality the optimal quantization 
yields a new optimization problem
whose optimal value is a lower bound for the optimal value of the original
problem (see \cite{lohndorfmodeling} for details).
Then, the exact lower bound
given by SDDP with resampling remains a lower bound for the exact lower
bound given by SDDP without resampling, which itself is by construction a lower bound
for the original problem.
In the numerical application, we fix the resampling size to $100$.
We denote by $\na{\VSDDP_t}_{t\in\ic{0,T}}$
the value functions returned by the SDDP algorithm.
Notice that, whereas the SDDP algorithm suffers from
the cardinality of the global uncertainty support, the DADP and PADP
algorithms do not.

We stop SDDP when the gap between its exact lower bound and a statistical
upper bound is lower than 1\%. That corresponds to the standard SDDP's
stopping criterion described in~\cite{shapiro10}, which is reputed
to be more consistent than the first stopping criterion introduced
in \cite{pereira1991multi}.
Our implementation of SDDP uses a level-one cut selection algorithm~\citep{guigues2017dual}
and keeps only the 100 most relevant cuts. By doing so, we significantly
reduce the computation time of SDDP.

\subsubsection{With the DADP and PADP algorithms}
to optimize the price and resource processes, we use a quasi-Newton method, and more precisely the L-BFGS algorithm
(implemented in the nonlinear solver Ipopt 3.12~\citep{wachter2006implementation}
which allows to explicitly tackle the linear
constraint $r\in \imag{(\INCIDENCEMATRIX)}$ in the resource decomposition
algorithm~\eqref{eq:district:subgradientquantdecom}).
The algorithms stop at~$\price$ and $\alloc$ either
when a stopping criterion is fulfilled or when no descent direction
is found.

\subsection{Devising online control policies}
\label{ssec:nodalpolicies}

Each algorithm (DADP, PADP and SDDP) returns a sequence of value
functions indexed by time. Using these value functions, we define
a sequence of \emph{surrogate global value functions}
$\na{\widehat{V}_t}_{t\in \ic{0,\final}}$ by\footnote{%
where the functions~\( \pvaluefunc_{\transport, t}\nc{\price} \)
(resp.~\( \qvaluefunc_{\transport, t}\nc{\alloc} \))
are easily deduced from~\eqref{eq:district:localarcpricepb}
(resp. from~\eqref{eq:district:localarcallocpb}).}
\begin{itemize}
\item
  $\widehat{V}_t =
  \VSDDP_t$ \hspace{0.0cm} for SDDP,
\item
  $\widehat{V}_t =
  \sum_{\node \in \NODES} \pvaluefunc_{\produ,t}^{\node}\nc{{\price}}+
  \pvaluefunc_{\transport, t}\nc{\price}$
  for DADP,
\item
  $\widehat{V}_t =
  \sum_{\node \in \NODES} \qvaluefunc_{\produ,t}^{\node}\nc{{\alloc}}+
  \qvaluefunc_{\transport, t}\nc{\alloc}$
  for PADP.
\end{itemize}
With these global value functions, we design online control policies.
For any time $t\in \ic{0,\final-1}$, any global state~$x_t \in \XX_t$
and global uncertainty $w_\post \in \WW_\post$,
the control policy is a solution of the following one-step optimization problem:
\begin{subequations}
  \label{eq:district:admissiblepolicy}
  \begin{align}
    \feedback_t
    & (x_t, w_\post) \in \argmin_{u_t}
      \Big(
      \min_{f_t, q_t}
      \nonumber
    \\
    & 
    \sum_{\node \in \NODES} L_t^{\node}(x_t^{\node}, u_t^{\node}, w^{\node}_\post) \! + \!
    \sum_{\arc \in \ARCS} l_t^{\arc}(q_t^{\arc}) \! + \!
      \widehat{V}_\post\np{x_\post} \Big)
      \nonumber
    \\
    &\st\ \quad \IncidenceMatrix q_t + f_t = 0 \eqfinv
    \\
    & \phantom{ \st\ \quad} x_\post^{\node} = \dynamic_t^{\node}(x_t^{\node}, u_t^{\node}, w^{\node}_\post) \eqfinv
    \\
    & \phantom{ \st\ \quad} \Delta_t^{\node}(x_t^{\node}, u_t^{\node}, w_\post^{\node}) = f_t^{\node} \eqsepv
    \forall \node \in \NODES \eqfinp
  \end{align}
\end{subequations}
As the policy induced by~\eqref{eq:district:admissiblepolicy}
is admissible for the global problem~\eqref{transportproblem},
the expected value of its associated cost is an upper bound
of the optimal value~$V_0\opt(x_0)$ of the original minimization
problem~\eqref{transportproblem}.

\subsection{Numerical results}
\label{Numerical_results}

We first compare the three algorithms \wrt\
the convergence and the CPU time needed for
computing offline value functions.
Second, we compare the values of the theoretical bounds for the optimal expected
total cost.
Third, we compare online policies simulation results.

\subsubsection{Computing offline value functions}
\begin{figure}[!ht]
  \centering
  \includegraphics[width=6.0cm]{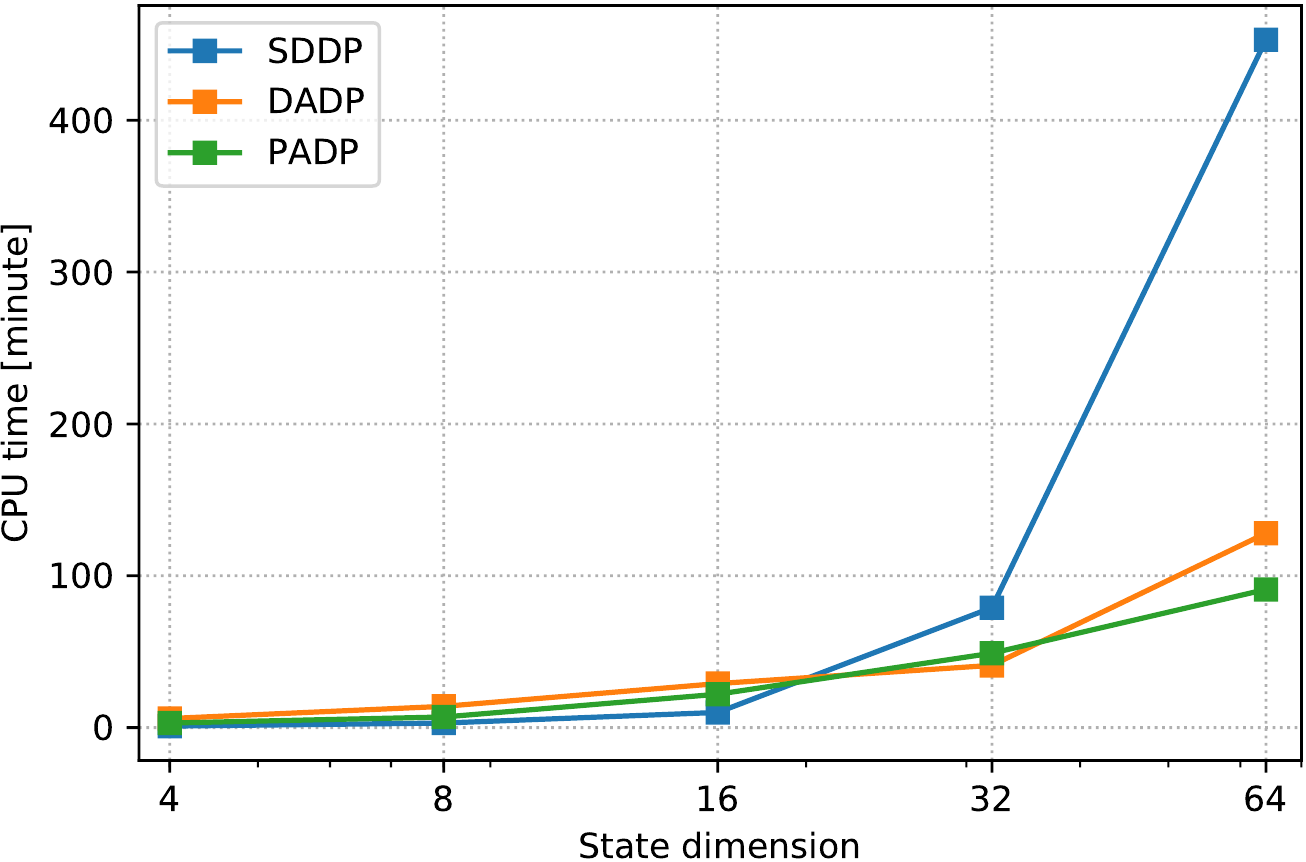}
  \caption{(a) CPU time for the three algorithms as a function
    of the state dimension }
  \label{fig:nodal:cputime}
\end{figure}

we solve Problem~\eqref{transportproblem} by SDDP,
price decomposition (DADP) and resource decomposition (PADP).
Table~\ref{tab:district:numeric:optres} details the execution
time and number of iterations taken before reaching convergence.
\begin{table}[!ht]
  \centering
  \resizebox{.35\textwidth}{!}{
    \begin{tabular}{|l|ccccc|}
      \hline
      \# Nodes $\card{\NODES}$    & \textrm{3}  \hspace{-0.2cm}
      & \textrm{6}  \hspace{-0.2cm}
      & \textrm{12} \hspace{-0.2cm}
      & \textrm{24} \hspace{-0.2cm}
      & \textrm{48} \hspace{-0.2cm} \\
      \hline
      dim~$\XX_t$          & 4                & 8                & 16
      & 32               & 64     \\
      \hline
      \hline
      SDDP CPU time      & 1'               & 3'               & 10'
      & 79'              & 453'   \\
      SDDP iterations    & 30               & 100              & 180
      & 500              & 1500   \\
      \hline
      \hline
      DADP CPU time      & 6'               & 14'              & 29'
      & 41'              & 128'   \\
      DADP iterations    & 27               & 34               & 30
      & 19               & 29     \\
      \hline
      \hline
      PADP CPU time      & 3'               & 7'               & 22'
      & 49'              & 91'    \\
      PADP iterations    & 11               & 12               & 20
      & 19               & 20     \\
      \hline
    \end{tabular}
  }
  \caption{Convergence results for SDDP, DADP and PADP}
  \label{tab:district:numeric:optres}
\end{table}
For a small-scale problem like \textrm{3-Nodes} (second column
of Table~\ref{tab:district:numeric:optres}), SDDP is faster
than DADP and PADP. However, for the 48-Nodes problem (last
column of Table~\ref{tab:district:numeric:optres}),
\emph{DADP and PADP} are \emph{more than three times faster}
than SDDP. Figure~\ref{fig:nodal:cputime} depicts how much CPU
time take the different algorithms with respect to the number
of state variables of the district. For this case study, we observe
that the \emph{CPU time grows almost linearly} \wrt\ the number
of nodes for DADP and PADP, whereas it grows exponentially for SDDP.
Otherwise stated, decomposition methods scale better than SDDP
in terms of CPU time for large microgrids instances.

\paragraph{Convergence of the SDDP algorithm}
on all instances, the approximate upper bound is estimated
every 10 iterations, with 1,000 scenarios. On the 12-Nodes problem, we observe that the gap
between the upper and lower bounds is below 1\% after 180 iterations and
that the lower bound remains stable after 250 iterations.

\begin{figure}[!ht]
  \begin{center}
    \begin{subfigure}[t]{\graphicscale\textwidth}
      {\includegraphics[width=\textwidth]{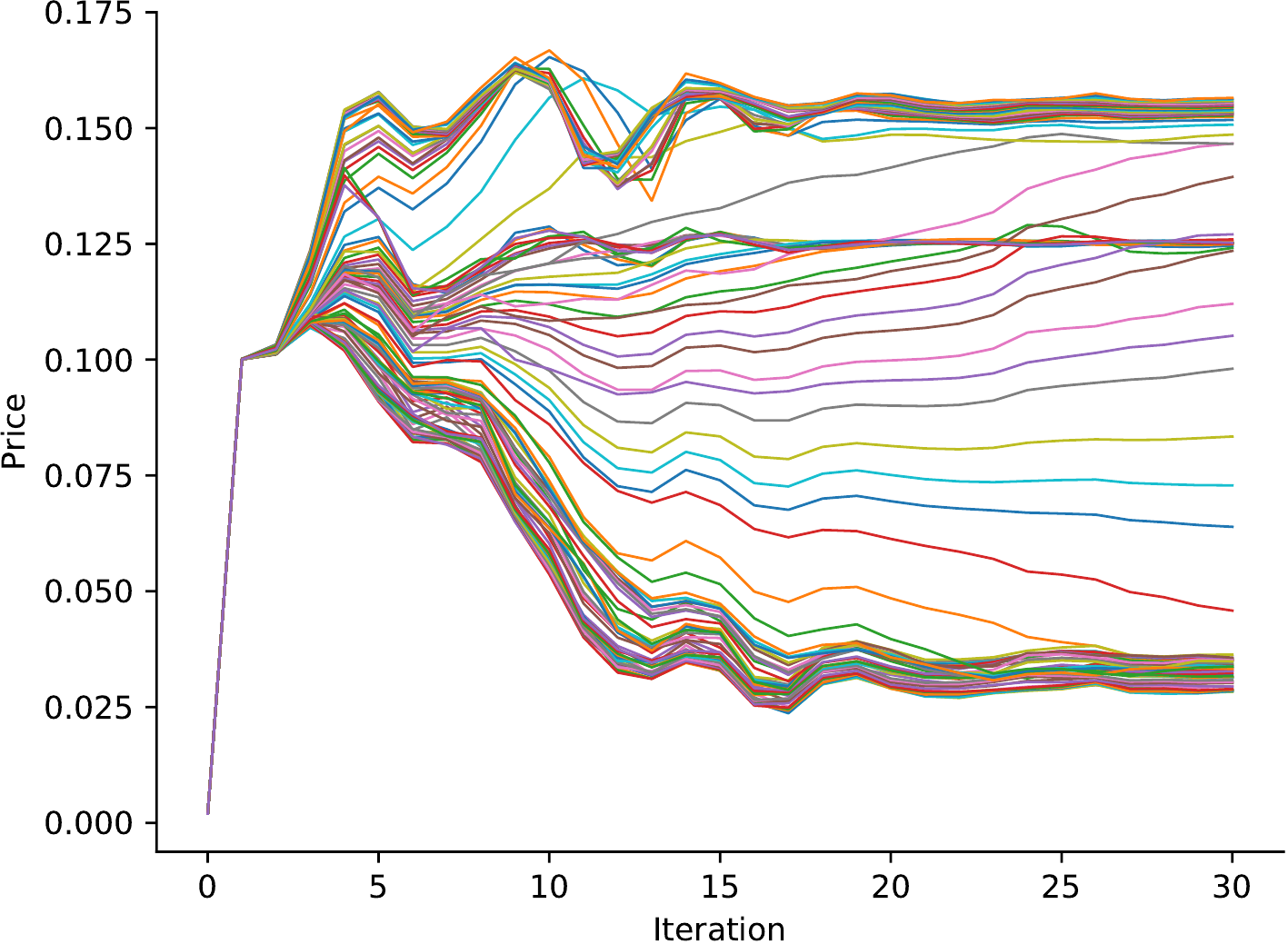}}
      \caption{\label{fig:subfigone}}
    \end{subfigure}\hspace{2mm}
    \begin{subfigure}[t]{\graphicscale\textwidth}
      {\includegraphics[width=\textwidth]{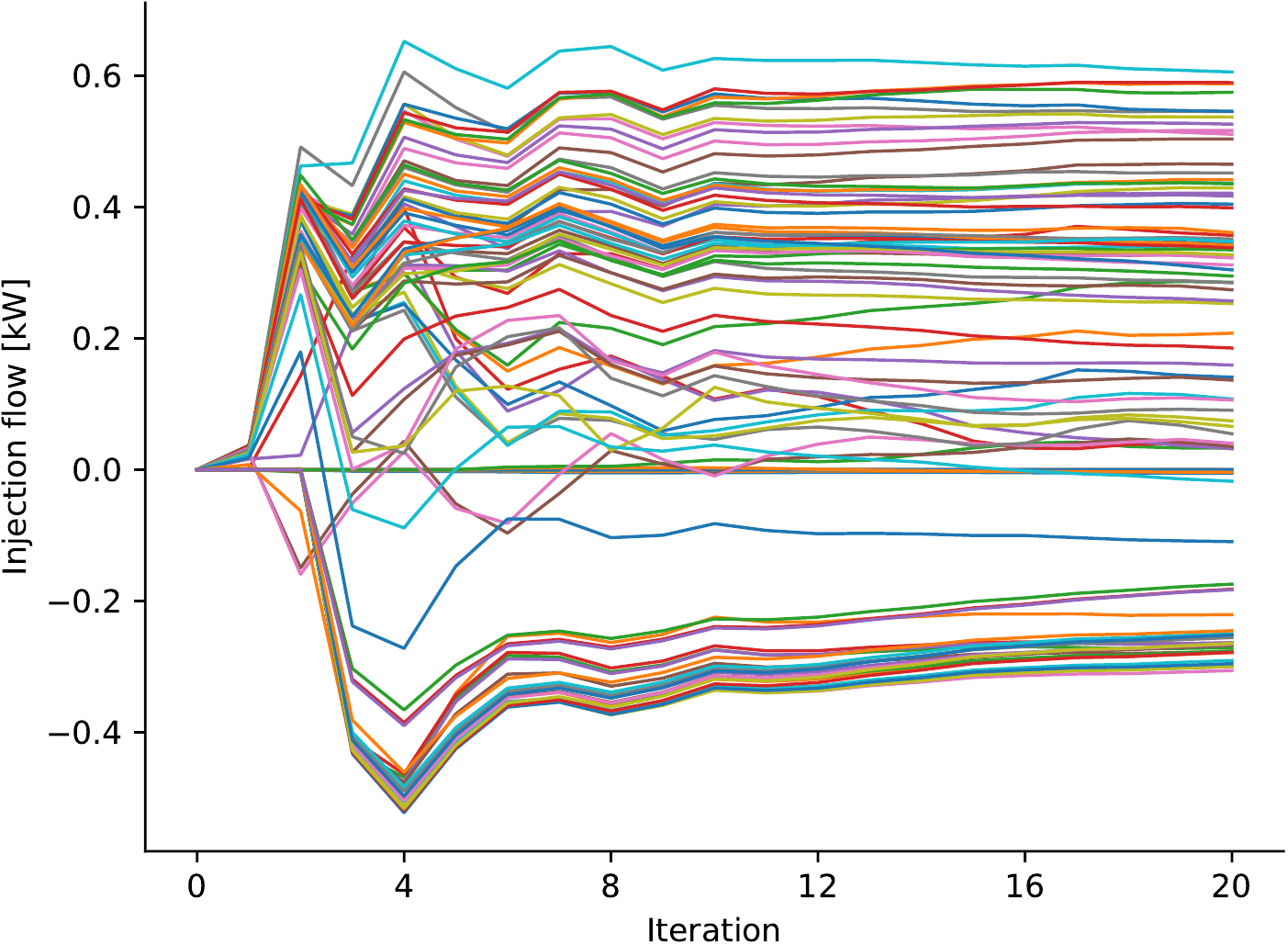}}
      \caption{\label{fig:subfigtwo}}
    \end{subfigure}
    \caption{Convergence of DADP prices (a) and PADP resources (b)
    for the 12-Nodes problem\label{fig:numerics:convdadp}}
  \end{center}
\end{figure}

\paragraph{Convergence of the DADP and PADP algorithms}
Figure~\ref{fig:numerics:convdadp}
shows the evolution of DADP's price process and PADP's resource process
over the iterations for the 12-Nodes problem.
We depict the convergence only
for the first node, the evolution of price process and resource
process in other nodes being similar.
On the left side of the figure, we plot the evolution
of the 96 different values of the price process
$\price^1 = (\price^1_0, \cdots,  \price^1_{\final-1})$
during the iterations of DADP. We observe that most of the prices
start to stabilize after 15 iterations, and do not exhibit significant
variation after 20 iterations. On the right side of the figure, we
plot the evolution of the 96 different values of the resource process
$\alloc^1 = (\alloc_0^1, \cdots, \alloc_{\final-1}^1)$ during the
iterations of PADP. We observe that the convergence of resources
is quicker than for prices, as the evolution of most resources
starts to stabilize after only 10 iterations.

\subsubsection{Theoretical bounds for the optimal expected total cost}
we then give the lower and upper bounds obtained by SDDP,
DADP, PADP in Table~\ref{tab:district:numeric:upperlower}.
The lower bound of the SDDP algorithm is the value
$\underline{V}_0^{sddp}(x_0)$ given by the SDDP method.
We recall that SDDP returns a lower bound because it uses
a suitable resampling of the global uncertainty distribution
instead of the original distribution itself (see the discussion
in~\S\ref{SDDP_on_the_global_problem}). DADP and PADP lower
and upper bounds are given by~\eqref{eq:district:relaxedmpts}
and~\eqref{eq:district:quantdecompositionrelaxed}
respectively. In Table~\ref{tab:district:numeric:upperlower},
we observe that SDDP and DADP lower bounds are close to each other,
and for problems with more than 12 nodes, DADP's lower
bound is up to 2.6\% better than SDDP's lower bound.
However, the gap between the upper bound given by PADP and
the two lower bounds is rather large.

\begin{table}[!ht]
  \centering
  \resizebox{.38\textwidth}{!}{
    \begin{tabular}{|l|ccccc|}
      \hline
      \# Nodes $\card{\NODES}$ & \textrm{3} & \textrm{6} & \textrm{12} & \textrm{24} & \textrm{48} \\
      \hline
      \hline
      SDDP LB  & 225.2      & 455.9      & 889.7       & 1752.8      & 3310.3      \\
      \hline
      DADP LB  & 213.7      & 447.3      & 896.7       & 1787.0      & 3396.4      \\
      \hline
      PADP UB  & 252.1      & 528.5      & 1052.3      & 2100.7      & 4016.6      \\
      \hline
    \end{tabular}
  }
  \caption{Upper and lower bounds, given by SDDP, DADP and PADP, for the optimal expected total cost}
  \label{tab:district:numeric:upperlower}
\end{table}
To sum up, DADP provides slightly better lower bounds than SDDP,
while being less computationally demanding
(and a parallel version of DADP would yield even
better performances).

\subsubsection{Online policies simulation results}
we now compare the performances of the different algorithms
in simulation. As explained in \S\ref{ssec:nodalpolicies},
we are able to devise online policies induced by SDDP, DADP
and PADP for the global problem, and to compute by Monte Carlo
an approximation of the expected cost of each of these policies.

The results obtained in simulation are given
in Table~\ref{tab:district:numeric:simulation}.
SDDP, DADP and PADP values are obtained by simulating the
corresponding policies on $5,000$ scenarios. The notation
$\pm$ corresponds to the 95\% confidence interval. We use
the value obtained by the SDDP policy as a reference,
and compute the relative difference:
a positive percentage means that the associated decomposition-based
policy is better than the SDDP policy.
Note that all these values correspond to admissible policies
for the global problem~\eqref{transportproblem}, and thus are
\emph{statistical} upper bounds of the optimal cost~$V_0\opt(x_0)$
of Problem~\eqref{transportproblem}.

\begin{table}
  \centering
  \resizebox{.5\textwidth}{!}{
    \begin{tabular}{|l|ccccc|}
      \hline
      \# Nodes  $\card{\NODES}$ & \textrm{3} & \textrm{6} & \textrm{12} & \textrm{24} & \textrm{48} \\
      \hline
      \hline
      SDDP value   & 226 $\pm$ 0.6    & 471 $\pm$ 0.8    & 936 $\pm$ 1.1     & 1859 $\pm$ 1.6    & 3550 $\pm$ 2.3    \\
      \hline
      \hline
      DADP value   & 228 $\pm$ 0.6    & 464 $\pm$ 0.8    & 923 $\pm$ 1.2     & 1839 $\pm$ 1.6    & 3490 $\pm$ 2.3    \\
      DADP/SDDP    & - 0.8 \%         & + 1.5 \%         & +1.4\%            & +1.1\%            & +1.7\%    \\
      \hline
      \hline
      PADP value   & 229 $\pm$ 0.6    & 471 $\pm$ 0.8    & 931 $\pm$ 1.1     & 1856 $\pm$ 1.6    & 3508  $\pm$ 2.2   \\
      PADP/SDDP    & -1.3\%           & 0.0\%            & +0.5\%            & +0.2\%            & +1.2\%    \\
      \hline
    \end{tabular}
  }
  \caption{Simulation results for 
    SDDP, DADP and PADP induced policies}
  \label{tab:district:numeric:simulation}
\end{table}

We make the following observations.
i) For problems with more than 6 nodes,
both the DADP policy and the PADP policy beat the SDDP policy.
ii) The DADP policy gives better results than the PADP policy.
iii) Comparing with the last line of Table
\ref{tab:district:numeric:upperlower}, the statistical
upper bounds obtained by the three simulation policies
are much closer to SDDP's and DADP's lower bounds than PADP's
exact upper bound. By assuming that the resource coordination
process is deterministic in PADP, we impose constant importation
flows for every possible realization of the uncertainties,
thus penalizing heavily the PADP algorithm (see also the
interpretation of PADP in the case of a decentralized
information structure in \cite[\S3.3]{Carpentier-Chancelier-DeLara-Pacaud:2020}).

\section{Conclusion}

We have addressed the mathematical problem of optimal management of urban microgrids by using
two decomposition algorithms relying
on a deterministic price (resp. resource) coordination process:
DADP and PADP. Both algorithms work in a distributed manner and are fully parallelizable.
We have conducted numerical simulations on microgrids of different sizes and
topologies, with  up to 48 buildings. We have compared the two decomposition
algorithms with a state-of-the-art SDDP algorithm.
Numerical results have shown the effectiveness of DADP, that gives better
results than the reference SDDP algorithm for problems with
more than 12~nodes --- both in terms of theoretical bounds and of
economic performance induced by online policies.
On problems with up to 48~nodes (corresponding to 64~state variables), we
observed that the performances of DADP and PADP scale well as the number of nodes
grew.
Numerically, we observe that decomposition-coordination methods are less impacted
by the curse of dimensionality than SDDP, as:
i) decomposed subproblems
have small dimension (1 or 2) and can be solved in parallel;
ii) as the size of the problem grows, we empirically observe that the number of iterations
of decomposition methods grows slower than with SDDP.
Thus, algorithms that mix spatial and temporal decompositions appear to be a
promising tool to address large-scale microgrid optimal management problems.

\end{document}